\documentclass{gtart}

\def\ifplaintex{\expandafter\ifx\csname documentclass\endcsname\relax}

\def\gtp{{\mathsurround=0pt\it $\cal G\mskip-2mu$eometry \&\ 
$\cal T\!\!$opology $\cal P\!$ublications}}  

\def\recd{{\small Received:\qua\receiveddate\ifx\reviseddate\relax
\else\qquad Revised:\qua\reviseddate\fi\par}} 

\def\lognumber#1{\def\thelognumber{#1}}
\def\volumenumber#1{\def\thevolumenumber{#1}}
\def\volumeyear#1{\def\thevolumeyear{#1}}
\def\papernumber#1{\def\thepapernumber{#1}}
\def\pagenumbers#1#2{\def\startpage{#1}\def\finishpage{#2}}
\def\published#1{\def\publishdate{#1}}

\def\received#1{\def\receiveddate{#1}}
\def\revised#1{\def\reviseddate{#1}}
\def\accepted#1{\def\accepteddate{#1}}
\def\asciititle#1{\def\theasciititle{#1}}

\def\asciiaddress#1{\def\theasciiaddress{#1}}

\long\def\asciiabstract#1{\long\def\theasciiabstract{#1}}

\let\\\par\let\thelognumber\relax\let\thevolumenumber\relax
\let\thepapernumber\relax\let\thevolumeyear\relax\let\startpage\relax
\let\finishpage\relax\let\publishdate\relax\let\receiveddate\relax
\let\reviseddate\relax\let\accepteddate\relax\let\theasciititle\relax
\let\theasciiauthors\relax\let\theasciiaddress\relax
\let\theasciiabstract\relax

\let\theasciiemail\relax

\ifplaintex
\font\logobig=cmssbx10 scaled 3836
\font\logomed=cmssbx10 scaled 2557
\else
\font\logobig=cmssbx10 scaled 4200
\font\logomed=cmssbx10 scaled 2800
\fi

\long\def\makeagttitle{   
\count0=\startpage
\agt\hfill      
\hbox to 45truept{\vbox to 0pt{\vglue -13truept{\logomed A\kern -.37em{\logobig 
T}\kern -.38em G}\vss}\hss}
\break
{\small Volume \thevolumenumber\ (\thevolumeyear)
\startpage--\finishpage\nl
Published: \publishdate}

\vglue .25truein

{\parskip=0pt\leftskip 0pt plus
1fil\def\\{\par\smallskip}{\Large\bf\thetitle}\par\medskip} \vglue
0.05truein

{\parskip=0pt\leftskip 0pt plus 1fil\def\\{\par}{\sc\theauthors}
\par\medskip}%
 
\vglue 0.03truein 

{\small\leftskip 25truept\rightskip 25truept{\bf Abstract}\stdspace\theabstract

{\bf AMS Classification}\stdspace\theprimaryclass
\ifx\thesecondaryclass\relax\else; \thesecondaryclass\fi\par
{\bf Keywords}\stdspace \thekeywords\par}\vglue 7truept

}   

\ifplaintex
\font\phead=cmsl9 scaled 950
\font\pnum=cmbx10 scaled 913
\font\pfoot=cmsl9 scaled 950
\headline{\vbox to 0pt{\vskip -4.5mm\line{\small\phead\ifnum
\count0=\startpage ISSN 1472-2739 (on-line) 1472-2747 (printed)
\hfill {\pnum\folio}\else\ifodd\count0\def\\{ }%
\ifx\theshorttitle\relax\thetitle\else\theshorttitle\fi\hfill{\pnum\folio}
\else\def\\{ and }{\pnum\folio}\hfill\ifx\theshortauthors\relax\theauthors
\else\theshortauthors\fi\fi\fi}\vss}}
\footline{\vbox to 0pt{\vglue 0mm\line{\small\pfoot\ifnum\count0=\startpage
\copyright\ \gtp\hfill\else
\agt, Volume \thevolumenumber\ (\thevolumeyear)\hfill\fi}\vss}}
\else
\font=cmsl9 scaled 1050
\font\lnum=cmbx10 
\font=cmsl9 scaled 1050
\makeatletter
\def\@oddhead{{\small\lhead\ifnum\count0=\startpage ISSN 1472-2739 
(on-line) 1472-2747 (printed)\hfill {\lnum\number\count0}\else\ifodd\count0
\def\\{ }\ifx\theshorttitle\relax \thetitle \else\theshorttitle\fi\hfill
{\lnum\number\count0}\else\def\\{ and }{\lnum\number\count0}
\hfill\ifx\theshortauthors\relax 
\theauthors\else\theshortauthors\fi\fi\fi}}\def\@evenhead{\@oddhead}
\def\@oddfoot{\small\lfoot\ifnum\count0=\startpage\copyright\ \gtp\hfill\else
\agt, Volume \thevolumenumber\ (\thevolumeyear)\hfill\fi}
\def\@evenfoot{\@oddfoot}
\makeatother
\fi
\let\maketitlepage\makeagttitle
\let\makeshorttitle\maketitlepage
\let\maketitle\maketitlepage

\newwrite\gtoutfile
\long\gdef\makeheadfile{  
{\def\\{, }\def\s{ }
\immediate\openout\gtoutfile head.xxx
\immediate\write\gtoutfile{To: math@arxiv.org}
\immediate\write\gtoutfile{Subject: put OR rep NNNNN:ppppp}
\immediate\write\gtoutfile{--text follows this line--}
\immediate\write\gtoutfile{Proxy-for: \ifx\theasciiauthors\relax
\theauthors\else\theasciiauthors\fi\s<\ifx\theasciiemail\relax\theemail\else\theasciiemail\fi>}
\immediate\write\gtoutfile{\noexpand\\}
\immediate\write\gtoutfile{Authors: \ifx\theasciiauthors\relax
\theauthors\else\theasciiauthors\fi}
{\def\\{ }\immediate\write\gtoutfile{Title: \ifx\theasciititle\relax
\thetitle\else\theasciititle\fi}}
\immediate\write\gtoutfile{Subj-class: GT or SG, GR etc}
\immediate\write\gtoutfile{MSC-class: \theprimaryclass\ifx\thesecondaryclass\relax\else, \thesecondaryclass\fi}
\immediate\write\gtoutfile{Journal-ref: Algebr. Geom. Topol. \thevolumenumber\s
(\thevolumeyear) \startpage-\finishpage}
\immediate\write\gtoutfile{Comments: Published by Algebraic and
Geometric Topology at}
\immediate\write\gtoutfile{\s\s\s  http://www.maths.warwick.ac.uk/agt/AGTVol\thevolumenumber/agt-\thevolumenumber-\thepapernumber.abs.html}
\immediate\write\gtoutfile{\noexpand\\}
\immediate\write\gtoutfile{}
\ifx\theasciiabstract\relax
\immediate\write\gtoutfile{\theabstract}\else
\immediate\write\gtoutfile{\theasciiabstract}\fi
\immediate\write\gtoutfile{}
\immediate\write\gtoutfile{\noexpand\\}
\immediate\write\gtoutfile{}
\immediate\closeout\gtoutfile}}  

\def\maketitlepage{\makeagttitle\makeheadfile}
\let\makeshorttitle\maketitlepage
\let\maketitle\maketitlepage

\lognumber{38}
\volumenumber{3}
\volumeyear{2003}
\papernumber{38}
\published{25 October 2003}
\pagenumbers{1089}{1101}
\received{12 September 2002}
\revised{7 July 2003}
\accepted{20 October 2003}

\usepackage{amsmath}
\usepackage{amssymb}
\usepackage{graphics}
\usepackage[all]{xy}

\theoremstyle{plain}
\newtheorem{theorem}{Theorem}[section]
\newtheorem{proposition}[theorem]{Proposition}
\newtheorem{lemma}[theorem]{Lemma}

\theoremstyle{definition}
\newtheorem{definition}[theorem]{Definition}
\newtheorem{example}[theorem]{Example}

\begin{document}

\title{Addendum to ``Coarse homology theories''}
\asciititle{Addendum to "Coarse homology theories"}

\author{Paul D. Mitchener}

\address{Institut f\"ur Mathematik, Universit{\"a}t 
G{\"o}ttingen\\D-37083 G\"ottingen, Germany} 
\asciiaddress{Institut fuer Mathematik, Universitaet 
Goettingen\\D-37083 Goettingen, Germany} 
\email{mitch@uni-math.gwdg.de}

\url{http://www.uni-math.gwdg.de/mitch/}

\begin{abstract}
This article corrects two mistakes in the article ``Coarse homology theories''
\cite{Mitch4}.
\end{abstract}

\asciiabstract{This article corrects two mistakes in the article
"Coarse homology theories" [Algebr. Geom. Topol. 1(2001) 271-297].}


\maketitle

\section{Introduction}

There are two mistakes in the article \cite{Mitch4}.  The first mistake is minor--- the definition of a coarsening cover is slightly too general for coarse homology theories to have the right properties.  Fortunately, this problem is easily fixed, and we can still prove an existence theorem concerning coarsening covers.

The second mistake is slightly more serious.  The original definition of a generalised ray is---as we show here---actually {\em too} general to be useful.  In this article we give a revised definition of a generalised ray that fixes this mistake; the basic philosophy of the earlier paper is still valid.  However, we are forced to amend our definition of a coarse $CW$-complex to be compatible with the new definition of a generalised ray.

In order to keep this paper relatively short we will not restate too many of the basic definitions from \cite{Mitch4}.  In particular, we assume that the reader knows what coarse spaces, coarse topological spaces, and coarse maps are.  We use without comment the product and disjoint union of coarse spaces defined in \cite{Mitch4}, as well as the quotient of a coarse space by an equivalence relation.

\subsection*{Acknowledgements}

The author wishes to thank Bernd Grave and Thomas Schick for valuable discussions.

\section{Generalised rays}

In \cite{Mitch4} an attempt was made to generalise the definition of the metric space $[0,\infty )$ in the coarse category.  In that paper, a generalised ray was defined to be the space $[0,\infty )$ equipped with some coarse structure compatible with the topology.  It is asserted in the final section of \cite{Mitch4} that any generalised ray, $R$, has trivial coarse $K$-homology (with coefficients in a $C^\ast$-algebra $A$), $KX_n (R;A)$, and that the $K$-theory of the coarse $C^\ast$-algebra, $C^\star_A (R)$, is trivial.  However, the following example shows that neither of these statements are true.

\begin{figure} 
\begin{center} 
\includegraphics{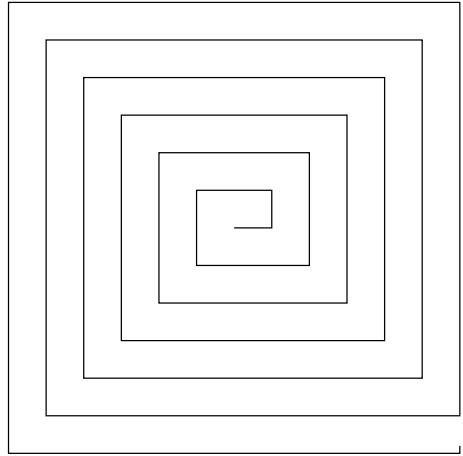} 
\end{center}
\caption{A ray embedded in ${\mathbb R}^2$} \label{spiral}
\end{figure}

\begin{example}
Let $R$ be the subset of ${\mathbb R}^2$ shown in figure \ref{spiral}.  Equip $R$ with the coarse structure inherited as a subset of the metric space ${\mathbb R}^2$.  Then $R$ is homeomorphic to the half-line $[0,\infty )$, and the given coarse structure is compatible with the topology.  However, it is clear that the spaces $R$ and ${\mathbb R}^2$ are coarsely equivalent.  Hence:
$$KX_n (R;A) = K_n (A) \qquad K_n C^\star_A (R) = K_n (A)$$
\end{example}

We therefore need a new, more restrictive, notion of a generalised ray.

\begin{definition}
Let $R$ be the space $[0,\infty )$ equipped with a unital coarse structure compatible with the topology.  We call the space $R$ a {\em generalised ray} if:

\begin{itemize}

\item Let $M,N\subseteq R\times R$ be entourages.  Then the set
$$M+N = \{ (u+x, v+y)\ |\ (u,v)\in M, (x,y)\in N \}$$
is an entourage.

\item Let $M\subseteq R\times R$ be an entourage.  Then the set
$$\overline{M} = \{ (u,v)\in R\times R\ |\ x\leq u,v\leq y,\ (x,y)\in M \}$$
is an entourage.

\item Let $M\subseteq R\times R$ be an entourage.  Then the set
$$\{ (x+a,y+a)\ |\ a\in R \}$$
is an entourage.

\end{itemize}

\end{definition}

Note that because the coarse structure is compatible with the topology, the subsets of a generalised ray that are {\em bounded} with respect to the coarse structure are precisely those that are bounded with respect to the metric.

\begin{proposition}
Let $R$ be a generalised ray.  Let $a\in R$.  Then the map $T_a \co R\to R$ defined by the formula $T_a (x) = x+a$ is close to the identity map $1_R$.
\end{proposition}

\begin{proof}
The map $P_a$ is coarse by definition of the ray.  By definition of a coarse structure, there is an entourage, $M\subseteq R\times R$, containing the point $(0,a)$.  By definition of a generalised ray, the set
$$\{ (x, T_a (x))\ |\ x\in R \}$$
is contained in entourage, which means, by definition, that the map $T_a$ is close to the identity map $1_R$.
\end{proof}

It follows that a coarse ray is flasque in the sense of \cite{HPR}.  In particular, the $K$-theory groups $K_n C^\star_A (X)$ are all trivial.

\begin{example}
We define the ray ${\mathbb R}_+$ be the space $[0,\infty )$ equipped with the coarse structure arising from the metric.  The entourages are subsets of {\em neighbourhoods of the diagonal}:
$$D_\alpha = \{ (x,y)\in {\mathbb R}_+ \ |\ |x-y|\leq \alpha \}$$
\end{example}

We shall reserve the notation ${\mathbb R}_+$ to denote the space $[0,\infty )$ equipped with the bounded coarse structure defined by the metric.

\begin{proposition} \label{M-S}
Let $R$ be a generalised ray.  Then every neighbourhood of the diagonal
$$D_\alpha = \{ (x,y)\in R\times R \ |\ |x-y| \leq \alpha \}$$
is an entourage.
\end{proposition}

\begin{proof}
Let $\Delta \subseteq R\times R$ denote the diagonal.  Then the set $\Delta \cup \{ (0,\alpha ) , (\alpha ,0 ) \}$ is an entourage.  By the second and third axioms in the definition of a generalised ray, the set
$$D_\alpha = \{ (x,y)\in R\times R \ |\ |x-y| \leq \alpha \}$$
must also be an entourage.
\end{proof}

It follows from the above proposition that the ray equipped with the $C_0$-coarse structure, as defined in \cite{Wr}, is not a generalised ray in our sense.

\begin{example}
Let $R$ be the space $[0, \infty)$.  Let $p_1 \co R\times R\to R$ and $p_2 \co R\times R\to R$ be the projections onto the first and second factors respecively.  Define a coarse structure by saying that an open subset $M\subseteq R\times R$ is an entourage if and only if for every point $x\in R$ the inverse images $p_1^{-1}(x)$ and $p_2^{-1}(x)$ are precompact (this coarse structure is in fact the continuously controlled coarse structure arising from the one point compactification of $R$).

The space $R$ is a generalised ray.
\end{example}

\begin{proposition}
The spaces ${\mathbb R}_+$ and $R$ are not coarsely equivalent.
\end{proposition}

\begin{proof}
The coarse structure on the space ${\mathbb R}_+$ is generated by a metric.  We will show that the space $R$ is not metrisable.

Suppose that the coarse structure on the space $R$ is generated by a metric, in the sense that there is a metric on $R$ such that every entourage is a subset of some uniformly bounded neighbourhood of the diagonal.  Then there is a sequence, $(M_n)$, of entourages such that every entourage $M\subseteq R\times R$ belongs to some member of the sequence $M_n$.

Choose points $(x_n , y_n )\in R\times R$ such that $(x_n , y_n )\not\in M_n$, $x_i\neq x_j$ for $i\neq j$, and $y_i \neq y_j$ for $i\neq j$.  Let
$$M = \bigcup_{n\in N} D((x_n,y_n)),1)$$  
where $D((x_n,y_n)),1)$  is the open disk of radius $1$ in the metric space $[0,\infty )\times [0,\infty )$ (say with the product metric).  Then according to the definition of the coarse structure on the space $R$, the open set $M$ is an entourage.  But there is no set in the sequence $M_n$ that contains $M$.

Therefore the coarse space $R$ is not metrisable, and we are done.
\end{proof}

\section{Coarse homology theories} \label{CHT}

Before we look at coarse homology theories, we should check exactly what we mean by coarse homotopy.  Actually, the notion is still essentially the same as that of \cite{Mitch4}, but we should be careful to use the new definition of a generalised ray.

\begin{definition}
Let $X$ and $Y$ be coarse spaces.  Let $f,g\co X\to Y$ be coarse maps.  Then a {\em coarse homotopy} linking $f$ and $g$ is a map $F\co X\times R\to Y$ for some generalised ray $R$ such that:

\begin{itemize}

\item The map $X\times R\to Y\times R$ defined by writing $(x,t)\mapsto (F(x,t),t)$ is a coarse map.

\item $F(x,0) = f(x)$ for every point $x\in X$.

\item For every bounded set $B\subseteq X$ there is a point $T\in R$ such that the function $F(x,t) = g(x)$ if $t\geq T$ and $x\in B$.

\item For every bounded set $B\subseteq X$ the set
$$\{ x\in X \ |\ F(x,t)\in B \textrm{ for some }t\in R \}$$
is bounded.

\end{itemize}

\end{definition}

The last condition in the definition of a coarse-homotopy did not appear in \cite{Mitch4}.  However, it is necessary for the homotopy-invariance arguments given in \cite{HPR, Mitch4} and earlier papers to work.  See \cite{Bart,MiS} for further discussion of this point.

More generally, we say that two coarse maps are {\em coarsely homotopic} if they are linked by a chain of coarse homotopies.

We now recall the main definition from \cite{Mitch4}.

\begin{definition} \label{cht}
A {\em coarse homology theory} consists of a collection of functors, $\{ HX_p \}_{p\in {\mathbb Z}}$, from the category of coarse spaces to the category of Abelian groups such that the following axioms hold:

\begin{itemize}

\item Coarse homotopy-invariance:

For any two coarsely homotopic maps $f\co X\to Y$ and $g\co X\to Y$, the induced maps $f_\ast \co HX_p (X)\to HX_p (Y)$ and $g_\ast \co HX_p (X)\to HX_p (Y)$ are equal.

\item Excision axiom:

Consider a decomposition $X = A\cup B$ of a coarse space $X$.  Suppose that for all entourages $m\subseteq X\times X$ we can find an entourage $M\subseteq X\times X$ such that $m(A)\cap m(B)\subseteq M(A\cap B)$.  Consider the inclusions $i\co A\cap B\hookrightarrow A$, $j\co A\cap B\hookrightarrow B$, $k\co A\hookrightarrow X$, and $l\co B\hookrightarrow X$.  Then we have a natural map $d \co HX_p(X)\to HX_{p-1}(A\cap B)$ and a long exact sequence:
$$\xymatrix@=8pt{
{} \ar[r] & HX_p (A\cap B) \ar[r]^-{\alpha} & HX_p (A)\oplus HX_p(B) \ar[r]^-{\beta} & HX_p (X) \ar[r]^-{d} & HX_{p-1}(A\cap B) \ar[r] & {}
}$$
where $\alpha = (i_\ast , -j_\ast )$ and $\beta = k_\ast + l_\ast$.

\end{itemize}

\end{definition}

A decomposition, $X = A\cup B$, of a coarse space $X$ is said to be {\em coarsely excisive} if the coarse excision axiom applies, that is to say for all entourages $m\subseteq X\times X$ we can find an entourage $M\subseteq X\times X$ such that $m(A)\cap m(B)\subseteq M(A\cap B)$.  The long exact sequence:
$$\xymatrix@=10pt{
{} \ar[r] & HX_p (A\cap B) \ar[r] & HX_p (A)\oplus HX_p(B) \ar[r] & HX_p (X) \ar[r] & HX_{p-1}(A\cap B) \ar[r] & {}
}$$
is called the {\em coarse Mayer-Vietoris sequence}.

The process of coarsening, described in \cite{HR2,Roe1}, is used to construct coarse homology theories on the category of proper metric spaces equipped with their bounded coarse structures.  This process can be generalised to more general coarse spaces as follows (see also \cite{STY}).

\begin{definition} \label{good}
Let $X$ be a coarse space.  A {\em good cover} of $X$ is a cover $\{ B_i \ |\ i\in I \}$ such that each set $B_i$ is bounded, and each set $B_i$ intersects only finitely many others in the cover.
\end{definition}

This differs slightly from the definition in \cite{Mitch4}.  For convenience, let us repeat definition 3.3 of \cite{Mitch4} where we are now using the above definition of good covers.

\begin{definition}
A directed family of good covers of $X$, $({\cal U}_i ,\phi_{ij} )_{i\in I}$, is said to be a {\em coarsening family} if there is a family of entourages $(M_i )$ such that:

\begin{itemize}

\item For all sets $U\in {\cal U}_i$ there is a point $x\in X$ such that $U\subseteq M_i (x)$.

\item Let $x\in X$ and suppose that $i<j$.  Then there is a set $U\in {\cal U}_j$ such that $M_i (x)\subseteq U$.

\item Let $M\subseteq X\times X$ be an entourage.  Then $M\subseteq M_i$ for some $i\in I$.

\end{itemize}

\end{definition}

The reason for our slight change of definition is that under the old definition of a good cover, proposition 3.6 of \cite{Mitch4} about the functoriality of coarse homology is actually incorrect.  However, everything is fine with the new definition.  To be precise, the following result is true.

\begin{theorem}
Let $\{ H^\mathrm{lf}_p \}$ be a generalised locally finite homology theory on the category of simplicial sets.  Then we can define a coarse homology theory on the category of coarse spaces that admit coarsening families by writing
$$HX_p (X) = \lim_{\to \atop i} H_p^\mathrm{lf} |{\cal U}_i |$$
where $X$ be a coarse space, with coarsening family $({\cal U}_i , \phi_{ij} )$.  
\qed
\end{theorem}

The proof of proposition 3.4 in \cite{Mitch4} about the existence of coarsening sequences is not valid with the above definition of a good cover.  However, we can prove a different existence result.

\begin{definition}
Let $X$ be a coarse space.  Then $X$ is said to have {\em bounded geometry} if it is coarsely equivalent to a space $Y$ where for every entourage $M\subseteq Y\times Y$, the number
$$\sup \{ |M(x)|\ |\ x\in Y \}$$
is finite.
\end{definition}

\begin{proposition}
Let $X$ be a coarse space of bounded geometry.  Then $X$ has a coarsening sequence.
\end{proposition}

\begin{proof}
Let us find a coarse space $Y$ equivalent to $X$ where for every entourage $M\subseteq Y\times Y$, the number
$$\sup \{ |M(x)|\ |\ x\in Y \}$$
is finite.  We will prove that the space $Y$ has a coarsening family.

Let $\{ M_i \ |\ i\in I \}$ be a cofinal family of entourages for $Y$ (in the sense that every entourage is contained in some entourage $M_i$) ordered by inclusion.  By hypothesis, we have a family of good covers, $\{ {\cal U}_i \ |\ i\in I \}$ defined by writing
$${\cal U}_i = \{ M_i (x) \ |\ x\in Y \}$$

But it is easy to check that this family is a coarsening family.
\end{proof}

\section{Coarse $CW$-complexes}

We begin by observing that the changed definition of a generalised ray means a small change in the definition of the building blocks of a coarse $CW$-complex.

\begin{definition}
Let $R$ be a generalised ray.  The {\em coarse $R$-sphere} of {\em dimension $n$} is the product $SX_R^n= (R\coprod R)^{n+1}$.  The {\em coarse $R$-cell} of {\em dimension $n+1$} is the product $DX_R^{n+1} = SX_R^n\times R$.  The coarse sphere
$$\{ (x,0)\ |\ x\in SX_R^n \}$$
is called the {\em boundary} of the coarse cell $DX_R^{n+1}$. 
\end{definition}

In particular, any generalised ray can be regarded as a coarse cell of dimension zero.  The disjoint union of two standard rays ${\mathbb R}_+$ is coarsely equivalent to the real line ${\mathbb R}$ with the bounded coarse structure coming from the metric.  If we think of a generalised ray as a `point at infinity', a disjoint union of two generalised rays appears as `two points at infinity'.  Generalising this idea to higher dimensions, we see that a coarse sphere is a `sphere at infinity' and a coarse cell is a `hemisphere at infinity'.

\begin{proposition} \label{Rhomotopy}
Let $R$ be a generalised ray.  Then the coarse map $i\co R\to (R\coprod R)^n\times R$ defined by the formula $i(s) = (0,s)$ is a coarse homotopy-equivalence.
\end{proposition}

\begin{proof}
Let $A\colon R\coprod R\rightarrow R$ be the coarse map that is equal to the identity map on each `copy' of the ray $R$ in the domain.  We then have a coarse map $p\colon (R\coprod R)^n\times R\rightarrow R$ defined by the formula
$$p(x_1 ,\ldots , x_n,s) = s+ \max (A(x_1), \ldots , A(x_n))$$

The composite $p\circ i$ is equal to the identity $1_R$.  Define a map $S\colon R\times R\rightarrow R$ by the formula $S(s,t) = \max (s-t,0)$, let $A = \max (A(x_1),\ldots , A(x_n))$, and write
$$H(x_1 , \ldots , x_n ,s,t) = \left\{ \begin{array}{ll}
(S(x_1,t) , \ldots , S(x_n,t),s+t) & t\leq A \\
(0, \ldots , 0 , p(x_1 , \ldots , x_n,s)) & t\geq A \\
\end{array} \right.$$

Then the map $H\colon (R\coprod R)^n\times R\times R\rightarrow R$ is a coarse homotopy between the composite $i\circ p$ and the identity $1_{(R\coprod R)^n\times R}$.
\end{proof}

Suppose we have a coarse space $Y$, and a coarse cell $DX^n$ with boundary $SX^{n-1}$.  If we have a coarse map $f\to SX^n\to Y$, we can form a new corse space $DX^n\cup_f Y$ by taking the quotient of the disjoint union $DX^n\coprod Y$ by the equivalence relation $x\sim f(x)$ for $x\in SX^n$.  The space $DX^n\cup_f Y$ is called the space $Y$ with an {\em attached} coarse cell.

\begin{definition}
A finite coarse $CW$-complex is a coarse space $X$ obtained by attaching a finite number of coarse cells to a finite disjoint union of generalised rays.
\end{definition}

It is clear that any finite coarse $CW$-complex has bounded geometry.

Let $\{ HX_p \}$ and $\{ HX'_p \}$ be coarse homology theories.  A {\em map} of coarse homology theories is a sequence of natural transformations $\alpha \co HX_n\to HX'_n$ that preserves coarse Mayer-Vietoris sequences.  We proved in \cite{Mitch4} that any map of coarse homology theories that is an isomorphism for generalised rays and one-point spaces is an isomorphism for finite coarse $CW$-complexes.  However, now that we have changed our definitions, we need to check that the argument of \cite{Mitch4} is still valid.

\begin{lemma} \label{sphere}
Let $\alpha \co HX_n (X)\to HX'_n (X)$ be a map of coarse homology theories  that is an isomorphism whenever the space $X$ is a generalised ray or the one point coarse space.  Then the map $\alpha$ is an isomorphism whenever the space $X$ is a coarse sphere.
\end{lemma}

\begin{proof}
We work by induction.  Let $R$ be a generalised ray.  Observe that the zero-dimensional sphere $SX_R^0$ is coarsely equivalent to a coarsely excisive union of generalised rays, $R_1\cup R_2$, and that the intersection $R_1\cap R_2$ is bounded and therefore equivalent to a single point, $+$.  We know that the maps $\alpha \co HX_n (+)\to HX'_n(+)$ and $\alpha \co HX_n (R_i)\to HX'_n (R_i)$ are isomorphisms.  An argument using Mayer-Vietoris sequences and the five lemma tells us that the map $\alpha \co HX_n (SX^0_R)\to HX'_n (SX^0_R)$ is an isomorphism.

Now, suppose that the map $\alpha \co HX_n (SX_R^{n-1})\to HX'_n (SX_R^{n-1})$ is an isomorphism.  We can write the coarse sphere $SX_R^n$ as a coarsely excisive union $D_1\cup D_2$, where $D_1$ and $D_2$ are coarse cells, and the intersection $D_1\cap D_2$ is coarsely equivalent to the sphere $SX_R^{n-1}$.  By proposition \ref{Rhomotopy} each cell $D_i$ is coarsely homotopy-equivalent to a generalised ray.  Therefore, by the same Mayer-Vietoris sequence argument as above, the map $\alpha \co HX_n (SX^n_R)\to HX'_n (SX^n_R)$ is an isomorphism.
\end{proof}

\begin{theorem} \label{CW}
Let $\alpha \co HX_n (X)\to HX'_n (X)$ be a map of coarse homology theories that is an isomorphism whenever the space $X$ is a generalised ray or the one point coarse space.  Then the map $\alpha$ is an isomorphism whenever the coarse space $X$ coarsely homotopy-equivalent to a finite coarse $CW$-complex.
\end{theorem}

\begin{proof}
The map $\alpha \co HX_n (X)\to HX'_n (X)$ is certainly an isomorphism whenever the space $X$ is a coarse $CW$-complex with just one cell.

Let $Y$ be a coarse $CW$-complex, and let $DX^n$ be a coarse cell with boundary $SX^n$.    Suppose that the map $\alpha \co HX_n (Y)\to HX'_n (Y)$ is an isomorphism, and we are given an attaching map $f\co SX^n\to Y$.  We must show that the map $\alpha \co HX_n (DX^n\cup_f Y)\to HX'_n (DX^n\cup_f Y)$ is an isomorphism.

Let $a\in [0,\infty )$.  The space $DX_a^n = \{ (x,t) \in DX^n \ |\ t\geq a \}$ is a coarse cell, and we have a coarsely excisive union:
$$DX^n\cup_f Y = (DX_1^n) \cup ((DX^n \backslash DX_2^n)\cup_f Y)$$ 

The space $(DX_2^n \backslash DX_1^n )\cup_f Y$ is coarsely equivalent to the space $Y$.  The intersection $(DX_1^n) \cap ((DX^n \backslash DX_2^n)\cup_f Y)$ is the space $DX_1^n \backslash DX_2^n$, which is coarsely equivalent to the coarse sphere $SX^{n-1}$.  Hence, by lemma \ref{sphere}, the desired result follows from an argument using Mayer-Vietoris sequences and the five lemma.  
\end{proof}

\section{The Novikov conjecture}

As we have already mentioned, there is a notion of a $C^\ast$-algebra, $C^\ast_A (X)$, associated to any coarse space $X$ and coefficient $C^\ast$-algebra $A$.  It is proved in \cite{HPR, Mitch4} that the sequence of functors $X\mapsto K_n C^\ast_A (X)$ is a coarse homology theory.

We have a locally finite generalised homology theory $X\mapsto KK^{-n}(C_0 (X),A)$ defined in terms of $KK$-theory.  We can coarsen it using the procedure described in section \ref{CHT} to define another coarse homology theory $X\mapsto KX_n (X;A)$ (at least when the space $X$ has bounded geometry).  There is a natural transformation of coarse homology theories
$$\alpha \co KX_n (X;A)\to K_n C^\star_A (X)$$
called the {\em coarse assembly map}.

\begin{lemma} \label{mainlemma}
Let $X$ be a topological space equipped with a proper continuous map $t\co X\to X$ such that:

\begin{itemize}

\item The map $t$ is properly homotopic to the identity map $1_X$.
 
\item For every compact subset $K\subseteq X$ there is a natural number $N$ such that $t^n [X]\cap K = \emptyset$ whenever $n\geq N$.

\item The family of induced maps $\{ t_\star^n \co C_0 (X)\to C_0 (X) \ |\ n\in {\mathbb N} \}$ is uniformly bounded.

\end{itemize}

Then the $KK$-theory groups $KK^{-n} (C_0(X),A)$ are all trivial.
\end{lemma}

\begin{proof}
We will prove that the $KK$-theory group $KK (C_0(X),A)$ is trivial for every $C^\ast$-algebra $A$.  The general result will then follow by Bott periodicity.  We naturally use an Eilenberg swindle.

Let $(H,F)$ be a Kasparov cycle for the pair $(C_0 (X),A)$.\footnote{See for example \cite{Bla} for more details concerning $KK$-theory.}  Thus $H$ is a Hilbert $A$-module equipped with a faithful representation of the $C^\ast$-algebra $C_0 (X)$ in the algebra of bounded linear operators ${\cal L}(H)$, and $F\in {\cal L}(H)$ is an operator such that the composites
$$(F^2 -1)\varphi \qquad (F-F^\star )\varphi \qquad F\varphi - \varphi F$$
are compact (in the sense of operators between Hilbert $A$-modules) for all functions $\varphi \in C_0 (X)$.

We have an induced map $t^\star \co C_0 (X)\to C_0 (X)$ such that the family $\{ (t^\star )^n \ |\ n\in {\mathbb N} \}$ is uniformly bounded, and for any given compactly supported function $\psi \in C_0 (X)$ the composite $(t^\star )^n (\varphi ) \psi$ is zero for all sufficiently large $n$, and all $\psi \in C_0 (X)$.  By the Hahn-Banach theorem the map $t^\star$ extends to a linear map $T \co {\cal L}(H)\to {\cal L}(H)$ such that the family $\{ T^n \ |\ n\in {\mathbb N} \}$ is uniformly bounded, and for any given compactly supported function $\psi \in C_0 (X)$ the composite $T^n (F) \psi$ is zero for all sufficiently large $n$.  Further, we can assume that the operator $T(F)$ is homotopic to the operator $F$.

We thus have a bounded operator
$$F^\infty = F\oplus T(F)\oplus T^2(F)\oplus \cdots$$
on the Hilbert space
$$H^\infty = H\oplus H\oplus H\oplus \cdots$$
If $\psi \in C_0 (X)$ is a compactly supported function, then all but finitely many terms in the series
$$F\psi \oplus T(F)\psi \oplus T^2(F)\psi \oplus \cdots$$
are zero.  It is now easy to verify that the pair $(H^\infty , F^\infty )$ is a Kasparov cycle.

Let $[(H,F)]$ be the element of the group $KK (C_0(X),A)$ represented by the cycle $(H,F)$.  Then certainly $[(H^\infty , T(F^\infty ))] = [(H^\infty , F^\infty )]$ and
$$[(H,F)] + [(H^\infty , T(F^\infty ))] = [(H^\infty , F^\infty )]$$ 

Therefore $[(H,F)]=0$ and we are done.
\end{proof}

\begin{theorem} 
Let $X$ be a coarse space coarsely homotopy-equivalent to a finite coarse $CW$-complex.   Then the coarse assembly map
$$\alpha \co KX_n (X;A)\to K_n C^\star_A (X)$$
is an isomorphism.
\end{theorem}

\begin{proof}
In view of theorem \ref{CW} we only need to prove the result when $X$ is a single point or a generalised ray.  The proof given in the paper \cite{Mitch4} when the space $X$ is a single point is fine.  However, the proof given of this result in \cite{Mitch4} when the space $X$ is a generalised ray is incorrect; we fix this mistake here.

Let $R$ be a generalised ray.  According to \cite{HPR} the $K$-theory groups $K_n C^\ast_A (R)$ are all zero since the space $R$ must be flasque.  We therefore need to prove that the groups $KX_n (R;A)$ are all zero. 

Let $S$ be the set of natural numbers, $\mathbb N$ (including zero) equipped with the coarse structure inherited as a subset of the ray $R$.  Then the spaces $S$ and $R$ are coarsely equivalent, and by proposition \ref{M-S} we can find an entourage, $M$, containing the set
$$\{ (i,j) \in {\mathbb N}\times {\mathbb N} \ |\ |i-j|\leq 1 \}$$

Let $(M_i)$ be an increasing family of entourages for the space $S$, ordered by inclusion, such that each entourage $M_i$ contains the enotourage $M$, and every entourage is contained in some entourage of the form $M_i$.  Define
$${\mathcal U}_i = \{ M_i (n) \ |\ n\in {\mathbb N} \}$$

Then the family of good covers, $( {\mathcal U}_i )$, is a coarsening family.  

Let $|{\mathcal U}_i |$ be the geometric realisation of the nerve of the cover ${\mathcal U}_i$.  Then the coarse $K$-homology group $KX_n (R;A)$ is by definition the direct limit of the groups $KK^{-n} (C_0 (|{\mathcal U}_i|),A)$.

Write $t (M_i (n) = M_i (n+1)$ for each vertex $M_i (n)$.  Then $t$ is a map from the set of vertices of the simplicial complex $|{\mathcal U}_i|$ to itself.  It can be linearly extended to a map $t\co |{\mathcal U}_i |\to |{\mathcal U}_i |$ by averaging over each simplex.

The following facts are clear from the definition of a generalised ray.

\begin{itemize}

\item The map $t$ is properly homotopic to the identity map $1_X$.
 
\item For every compact subset $K\subseteq X$ there is a natural number $N$ such that $t^n [X]\cap K = \emptyset$ whenever $n\geq N$.

\item The family of induced maps $\{ t_\star^n \co C_0 (X)\to C_0 (X) \ |\ n\in {\mathbb N} \}$ is uniformly bounded.

\end{itemize}

Therefore, by lemma \ref{mainlemma} we are done.
\end{proof}

The applications of the above theorem to the Novikov conjecture described in \cite{Mitch4} still work, although we need the new definition of a coarse $CW$-complex featuring in this article.

\Addresses\recd

\end{document}